\documentclass[12pt]{article}
\usepackage{amsmath,amsthm,amsfonts,amssymb,amscd}
\usepackage[latin1]{inputenc}
\usepackage[brazil]{babel}

\newtheorem{theorem}{Theorem}[section]
\newtheorem{coro}{Corollary}

\newtheorem{proposition}{Proposition}

\begin{document}


\begin{center}
{\Large Classes of Weingarten  Surfaces in ${\mathbb{S}}^2\times {\mathbb{R}}$}
\end{center}

\begin{center}
{Armando  V. Corro, Marcelo  A. Souza, Romildo Pina}
\end{center}
\footnote{The authors' research were partially supported by FAPEG and CAPES/PROCAD - NF.}

{corro@ufg.br, msouza@ufg.br, romildo@ufg.br}

\begin{abstract}
{In this work we  study surfaces in radial conformally flat 3-spaces. 
We characterize surfaces of rotation with constant Gaussian and Extrinsic curvature in these radial 3-spaces. 
We prove that all the spheres in the conformal 3-space have constant Gaussian curvature  $K=1$  if, and only if, the conformal factor is special. In this special case we study geometric properties of this ambient 3-space, and as an application we prove that it is isometric to the space ${\mathbb{S}}^2\times {\mathbb{R}}$, so we consider it as the {\em Radial Model} of ${\mathbb{S}}^2\times {\mathbb{R}}$.   
We obtain  two  classes of Weingarten surfaces in the {\em Radial Model}, which satisfy $\tilde{K}_E+\tilde{H}^2-\tilde{K}=0 $ and $2\tilde{K}_E-\tilde{K}=0 $, where $\tilde{K}$ is the Gaussian curvature, $\tilde{H}$ is the mean curvature and $\tilde{K}_E$ is the extrinsic curvature. Moreover, by using the relations between the curvatures of the {\em Radial Model} and the curvatures with respect to the euclidean metric ([CPS]), we  prove that  first class the Weingarten  surfaces in {\em Radial Model} corresponds, up  to isometries, to the minimal surfaces in $\mathbb{R}^3$, and second class corresponds to EDSGHW - surfaces in Euclidean space $ \mathbb{R} ^ 3$(\cite{DC}). Consequently these two classes of surfaces  have a Weierstrass type representation depending  on two holomorphic functions.
}
\end{abstract}
{keyword:}
{radial conformal metrics;  extrinsic curvature; Weingarten surfaces}


\section{Introduction}
Surfaces with constant Gaussian curvature have been studied by several people (see \cite{AEG}, \cite{BKP}, \cite{CMM},
\cite{CPR}, \cite{MO}). A lot of recent works on product spaces have meant a renewed interest in these
3-manifolds and have turned this research topic into a fashionable one. An up to date reference list of
papers on this subject can be found in \cite{FM}. Recently the study of surfaces  with constant extrinsic
curvature has been  extensively incremented, for example, the paper   \cite{EGR} for surfaces in product spaces.
In \cite{EGR} the authors proved that every complete connected immersed surface with positive extrinsic
curvature in ${\mathbb{H}}^2\times {\mathbb{R}}$ must be properly embedded,  and  homeomorphic to a sphere or a plane, respectively,
and, in the latter case, they studied the behavior of the end. They focused  their attention on surfaces with
constant positive  extrinsic curvature. The authors \cite{CPS} proved the existence of complete surfaces of rotation 
with non positive constant extrinsic curvature in a conformally flat 3-space. 
In  \cite{AEG}, the authors studied surfaces with constant Gaussian
curvature in ${\mathbb{S}}^2\times {{\mathbb{R}}}$ and ${\mathbb{H}}^2\times {{\mathbb{R}}}$.

The weingarten surfaces has been extensively studied in the present day. In \cite{FY} the authors  studied helicoidal Weingarten surfaces in the 3-dimensional Euclidean space. In \cite{MF}  the authors  classify  the complete rotational special Weingarten surfaces in  ${\mathbb{S}}^2\times {{\mathbb{R}}}$ and ${\mathbb{H}}^2\times {{\mathbb{R}}}$,  whose mean curvature $H$  and extrinsic curvature $K_e$ satisfy  $ H = f(H^2 - K_e)$, for some function $f$. Furthermore we show the existence
of non-complete examples of such surfaces. In \cite{MF}  the author provided a vertical height estimate for compact
special Weingarten surfaces of elliptic type in  ${\mathbb{M}}^2\times {{\mathbb{R}}}$, i.e. surfaces whose mean curvature $H$  and extrinsic Gauss curvature $K_e$  satisfy 
$ H = f(H^2 - K_e)$  with $ 4x{f^{'}(x))}^2  <  1$, $ \forall x \geq 0$.  

In \cite{DC}, the authors introduced the study of  classes  of generalized Weingarten surfaces  in the Euclidean space $\mathbb{R}^3$. In particular it was studied one class of surfaces which satisfy the relation  $ <X,X> K +2 <X,N>H =0$. These surfaces are called of EDSGHW - surfaces, moreover  they obtained for them a Weierstrass type representation depending  on two holomorphic  functions.

The study of surfaces in  spaces that are conformal to the Euclidean space is natural,
because they include the spaces of constant curvature, the punctured sphere ${\mathbb{S}}^3\setminus \{p\}$ and the hyperbolic space ${\mathbb H}^3$. Then it is natural to
consider surfaces in some special spaces with
a metric conformal  to the Euclidean metric. In this paper we consider the space ${\mathbb{R}}^3_g:=({\mathbb{R}}^3, g_F)$, where
$g_F:=(g_{ij})=(\delta_{ij})/F^2$ and $F=F(t)\neq 0$ is a differentiable function, where $t=x_1^2+x_2^2+x_3^2$,   $x= (x_{1},x_{2},x_{3})\in
{\mathbb{R}}^3$. We observe that if $F$ is bounded then the ambient space  ${\mathbb{R}}^3_g$ is a complete
  Riemannian manifold.

The above spaces were studied in  \cite{CPS}, in the particular case which $F(t)=exp(-t)$ the space ${\mathbb{R}}^3_g$ 
was denoted by ${\mathbb{E}}_3$.  This particular metric appears  as
a solution to the Einstein equation obtained by Pina and Tenenblat \cite{PT},
 with a great potential of applications in physics (see \cite{SKMHH} and the  references there).   
  As the space ${\mathbb{R}}^3_g$ is invariant under the actions of the
 orthogonal group,  it is natural to consider surfaces  of rotation that are invariant under the same group.

The main purpose of this  work is to study surfaces of rotation in ${\mathbb{R}}^3_g$. We show that the round spheres $S(0,R)$
of radius $R$, centered at origin, have   constant extrinsic curvature and that the radial lines are geodesics in ${\mathbb{R}}^3_g$.
We proved that the round spheres centered at origin $S(0,R)$
 have zero extrinsic curvature if, and only if, the conformal factor in the ambient space is given by
 $F(t)=\sqrt{t}$. 

In this special case we study geometric properties of the ambient 3-space, $({\mathbb{R}}^3\setminus \{0\}, g_F)$, $F(t)= \sqrt{t},$ it is a complete manifold with non negative sectional curvatures, and 
all the spheres $S(0,R)$ are compact minimal surfaces and have constant Gaussian curvature  $K=1$. As an application we prove that it is isometric to the space ${\mathbb{S}}^2\times {\mathbb{R}}$, so we consider it as the {\em Radial Model} of ${\mathbb{S}}^2\times {\mathbb{R}}$. We can use this fact  to study problems related to the classification of surfaces with special properties in this new model.
We obtain  two  classes of Weingarten surfaces in the {\em Radial Model}, which satisfy $\tilde{K}_E+\tilde{H}^2-\tilde{K}=0 $ and $2\tilde{K}_E-\tilde{K}=0 $, where $\tilde{K}$ is the Gaussian curvature, $\tilde{H}$ is the mean curvature and $\tilde{K}_E$ is the extrinsic curvature. Moreover, by using the relation between the curvatures of the {\em Radial Model} and the curvatures with respect to the euclidean metric ([CPS]), we  prove that  first class of the Weingarten  surfaces in {\em Radial Model} corresponds, up  to isometries, to the minimal surfaces in $\mathbb{R}^3$, and second class corresponds to EDSGHW - surfaces in Euclidean space $ \mathbb{R} ^ 3$. Consequently these two classes of surfaces  has a Weierstrass type representation depending  on two holomorphic functions.

\section{Preliminaries}

        We consider $({\mathbb{R}}^3, <\;,\;>)$, the 3-dimensional real vector space equipped with the canonical Euclidean metric, 
and  $({\mathbb{R}}^3, <\;,\;>_g)$, the Euclidean 3-dimen\-sional space equipped with a metric that is conformal to the Euclidean me\-tric. 

In the next result we will consider a regular parametrized surface $X=X(u,v)$, where $(u,v)$ are Isothermal parameters of $X$, namely  $E=<X_u,X_u>=<X_v,X_v>=G$ and $<X_u,X_v>=0$. 

\begin{proposition} { Let $X=X(u,v)$ be a regular parametrized surface and assume that $X$ is isothermal.
 Then: 

\begin{equation}
\label{eq9}
X_{uu}+X_{vv}=2EHN,
\end{equation}\\
The Gaussian curvature is given by
\begin{equation}
\label{eq1}
K=\frac{-1}{2E}\left\{\left( \frac{{E}_v}{{E}} \right)_v+\left( \frac{{E}_u}{{E}} \right)_u \right\},
\end{equation}
where  $H$ is the mean curvature and $N$ is the Gauss map. 
The components of the metric $<\;,\;>_g$ are given by
    \begin{equation}
    \label{gdeltaijF2}
    g_{ij}(x)={\delta_{ij}}/{F^2(x)}, \,\,\,\,x=(x_1,x_2,x_3),\,\, 1\leq i,j\leq 3,
    \end{equation}
    where $F: {\mathbb{R}}^3 \longrightarrow {\mathbb{R}}\setminus \{0\}$ is a differentiable function. }
		\end{proposition}
		
We observe that if $F$ is bounded then the conformal metric $<\;,\;>_g$ is a
complete metric. 

Considering  the  Levi Civita connection $\bar{\nabla}$ of  $({\mathbb{R}}^3,<\;,\;>_g)$
and the canonical basis  $\{e_1,e_2,e_3\}$ of
${\mathbb{R}}^3$, we get
\begin{equation}
\label{nablaBarra}
\bar{\nabla}_{e_i}{e_j}=\sum_{k=1}^{3} \Gamma_{ij}^k e_k,\,\,\mbox{and } \bar{\nabla}_{e_i}{e_j}=
\bar{\nabla}_{e_j}{e_i}.
\end{equation}

Since, in Equation (\ref{gdeltaijF2}),  $g_{ij}=0$ for $i\neq j$,  the Christoffel's
 symbols of this metric are given by:
   \begin{equation}
    \label{SimbChristoffel2}
    \begin{array}{lcl}
     \Gamma_{ij}^k&=&0,\,\,\,\, i\neq j\neq k\neq i,\\
     \Gamma_{ii}^j&=&{\displaystyle {F_{,j}}/{F}, \,\,\, \forall\,\, i\neq j,}\\
     \Gamma_{ij}^i&=&{\displaystyle {-F_{,j}}/{F},\,\,\, 1\leq i,j\leq 3,}
    \end{array}
    \end{equation}
where $F_{,j}={\partial F}/{\partial x_j}$ denotes the partial derivative of $F$ with respect to $x_j$.

The Riemannian manifold $({\mathbb{R}}^3,<\;,\;>_g)$  has sectional curvature given by (see \cite{dC})
 \begin{equation}
\label{KsecGeral}
 K\left({\partial}/{\partial x_i},{\partial}/{\partial x_j}\right)(x)= \left[\left({F_{,i}}/{F}\right)_{,i}+\left({F_{,j}}/{F}\right)_{,j} -\left({F_{,k}}/{F}\right)^2 \right]F^2,
  \end{equation}
 where $ 1 \leq i,k , j\leq 3$, are distincts.

Let $\beta(s)=(x_1(s),x_2(s),x_3(s))$ be a curve parametrized by arc length. The curve $\beta(s)$ is
a geodesic if, and only if, the component functions of $\beta$ satisfy the system of ordinary differential equation 
   \begin{equation}
    \label{Geodesicas01}
    \displaystyle \frac{d^2 x_k}{ds^2}+\sum_{i,j}\Gamma_{ij}^{k}\frac{dx_i}{ds}\frac{dx_j}{ds}=0, 
    \end{equation}
$k=1,2,3$, where the Christoffel's symbols are given by (\ref{SimbChristoffel2}).
    
We know that given a $X$ parametrized surface in the Euclidean 3-space, the mapping $X_I =\frac{X}{<X,X>}$ is the inversion of $X$ with respect to the origin,  with normal Gauss map $N_I $  given by 

$$  
N_I =\frac{-2<X,N>}{<X,X>}X + N 
$$

 and that the Weingarten matrices are related by
 \begin{equation}
  \label{matrizw}
    W_I= <X,X> W-<X,N>2I,
     \end{equation}

where $ W_I$  and $W$  are the Weingarten matrices of the $X_I$  and $X $ respectively.

  We note that the mean curvatures $H_I$  and $H$ of the  $X_I$  and $X$  respectively, satisfies:

 \begin{equation}
  \label{matrizwH}
<X,X>H + 2<X,N>=H_I
\end{equation}

The following relations were motivated by \cite{C}, where they were obtained  in the case that the ambient space is the
 hyperbolic space ${\mathbb H}^3$, the proofs can be found in \cite{CPS}.
\begin{theorem}
{\cite{CPS}} {  Let $X:U \subset {\mathbb{R}}^{2} \rightarrow
{\mathbb{R}}^{3}$ be a  regular parameterized surface. Consider $X(U)$
as a surface in $({\mathbb{R}}^{3},<\;,\;>)$
with the Euclidean metric,  let  $N$ be the normal Gauss mapping, $-\lambda_i$
the principal curvatures, $H$ and $K$ the  mean and Gaussian
curvatures, respectively. Analogously, consider $X(U)$ like a
surface in $({\mathbb{R}}^{3},<\;,\;>_g)$, with a  metric conformal to the
Euclidean metric, with the conformal factor $F^{-2}$, let
$-\tilde{\lambda_i}$ be the
principal curvatures, $\tilde{H}$ and $\tilde{K}_E$ the  mean and
the  extrinsic curvatures, respectively. Then
$$\begin{array}{lcl}
\tilde{\lambda_i}&=&F\lambda_i - <N,grad F>,
\\
\tilde{H}&=&FH+<N,grad F>,\\
\tilde{K}_E&=&F^2K+2HF<N,grad F>+<N,grad F>^2,
\end{array}
$$
where $F$ denotes the  evaluation of $F$ at $X(u,v), \;\;(u,v) \in U$.}
\end{theorem}

\section{Surfaces with constant curvatures in  conformally flat spaces}

In this section we consider radial conformal metrics. Our aim is study surfaces of rotation with constant extrinsic
curvature in   ${\mathbb{R}}^3_g$. 
\begin{proposition}
   {The Riemannian manifold ${\mathbb{R}}^3_g$ has sectional curvature given by }
 \begin{equation}
\label{Ksec}
 K\left({\partial}/{\partial x_i},{\partial}/{\partial x_j}\right)(x)= {-4 \dot{F}^2x^2_{k}}+{4F\dot{F}}+{4(x^2_i+x^2_j)(-\dot{F}^2+F\ddot{F})},
  \end{equation}
{  where $ 1 \leq i,k , j\leq 3$, are distincts, and if $\beta(s)=(x_1(s),x_2(s),x_3(s))$ is a curve parametrized by arc lenght,  then $\beta(s)$ is a geodesic if, and only if, }
   \begin{equation}
    \label{Geodesicas2}
    \displaystyle \frac{d^2 x_k}{ds^2}+\frac{2x_k\dot F}{F}\left[\sum_{i\neq k}\left(\frac{dx_i}{ds}\right)^2-\left(\frac{dx_k}{ds}\right)^2\right] -\frac{4\dot F}{F}\frac{dx_k}{ds}
    \sum_{i\neq k} x_i \frac{dx_i}{ds}=0, 
    \end{equation}
{ $k=1,2,3,$  where  $\dot{F}$ means the derivative of $F$ with  respect to the variable $t$.}
\end{proposition}
PROOF.
The proof follows by straightforward calculations.  The Equation (\ref{Ksec}) follows from Equations (\ref{SimbChristoffel2}) and (\ref{KsecGeral}), and from the fact that$ F,_{j}=\dot{F}(t)2x_j.$
To obtain the Equation (\ref{Geodesicas2}), simply replace   the Christoffel's symbols (\ref{SimbChristoffel2}) in equation (\ref{Geodesicas01}).
\hfill $ \Box $

We remark that the orthogonal maps are isometries of the space ${\mathbb{R}}^3_g$, then unless of isometries any surface of rotation around an axis through the origin we can, without loss of generality, consider the parametrization of a surface of rotation given by $X(u,v)=(\varphi(u)\cos(v),\varphi(u)\sin(v),u)$.
\begin{proposition}
 {Let $X(u,v)=(\varphi(u)\cos(v),\varphi(u)\sin(v),u)$ be a surface of rotation in ${\mathbb{R}}^3_g$. The surface $X$
has constant extrinsic curvature $c_0$ if, and only if, $\varphi$ satisfies the ordinary differential equation }
\begin{equation}
\label{eqC}[F+2\varphi \dot{F}(-\varphi+u\varphi')][F\varphi''-2a^2\dot{F}(-\varphi+u\varphi')]=-c_0a^4\varphi,
\end{equation}
{\it where} $a^2=1+\varphi'^2.$
\end{proposition}
PROOF.
  The coefficients of
 the first fundamental form, with respect to the Euclidean metric $ g_{0}$, are given by
$$ <X_u,X_u>_{g_0}=1+\varphi'^2(u),  \quad \quad  <X_v,X_v>_{g_0} =\varphi^2(u) \quad \mbox{ and } \quad
  <X_u,X_v>_{g_0}=0.$$
The coefficients of the second fundamental form are given by
$$e={-\varphi''(u)}/{a},\;\; g={\varphi(u)}/{a}, \;\; \mbox{ and } \;\; f=0.$$
We observe that
\begin{equation}
\label{XinternoN} <X,N> = {-\varphi+u\varphi'}/{a},\quad \quad \quad
<X,X> = \varphi^{2}+ u^{2}.
\end{equation}

\noindent Hence $\; <N,grad \,F> = 2\dot{F} < N,X>$. Then, using Theorem 2.1  and equation
(\ref{XinternoN}),  we have
$$\tilde{\lambda_i}=F(t)\lambda_i - 2<N,X> =
F(t)\lambda_i-2\dot{F}\left({-\varphi+u\varphi'}\right)/{a},$$
where $\lambda_{1}={\varphi''}/{a^3}$ e $\lambda_{2}=
-{1}/{\varphi a}$. In this case the extrinsic curvature
$\tilde{K}_E= \tilde{\lambda_{1}}\tilde{\lambda_{2}}$ is a constant
$c_{0}$ if, and only if,  equation (\ref{eqC}) is satisfied. This concludes the proof of Proposition 3.
  \hfill $ \Box $

In the next result, we present an explicit solution for the equation (\ref{eqC}), when 
$F(t)=exp(-t)$ and $c_0>0$. This special case were studied in \cite{CPS},
 where we have considered $c_0\leq 0$, and obtained explicit solutions for (\ref{eqC}).

\begin{proposition}
 {Let $X(u,v)=(\varphi(u)\cos(v),\varphi(u) \sin (v),u)$ be a surface of rotation in ${\mathbb{E}}^3$. We have that $X$ have constant extrinsic curvature $c_0$ if, and only if, $\varphi$ satisfies the ordinary differential equation}
\begin{equation}
\label{eqGaussCte}[1+2\varphi^2 -2u\varphi \varphi']\varphi''+a^2[4\varphi(\varphi')^2+2\varphi+2u\varphi']=-c_0
a^4\varphi e^{2
(u^2+\varphi^2)},
\end{equation}
{where $a^2=1+(\varphi')^2.$  The ordinary differential equation (\ref {eqGaussCte}) admits the explicit solution 
$\varphi(u) = \sqrt{R^2- u^2}$, where $R$ is uniquely determinated by the relation}

\begin{equation}
\label{Czero}
c_{0} = {\left(1+2R^2\right)^2}/{R^2 e^{2R^2}}. 
\end{equation}
\end{proposition}
PROOF.
To obtain equation (\ref{eqGaussCte}) we use equation (\ref{eqC}) and the fact that \linebreak $\dot{F} = -F=-e^{-({u^2+\varphi^2(u)})}$.
By straightforward calculation $\varphi(u) = \sqrt{R^2- u^2}$ is a solution of (\ref{eqGaussCte}), where $R$ and $c_0$ are related by equation (\ref{Czero}).
 \hfill $ \Box $

 We remark that the surface generated by rotating the curve that is graph of the function
 $\varphi(u)= \sqrt{R^2- u^2}$, around the $u-$axis, is a sphere centered
 at origin with radius $R$. 

Now, in the next theorem, we will  show that all the Euclidean spheres $S(0,R)$
have constant non negative extrinsic curvature in ${\mathbb{R}}^3_g$.
\begin{theorem}
 {The spheres $S(0,R)$ in ${\mathbb{R}}^3_g$ have  non negative constant extrinsic curvature given by
$$\tilde{K}_E= \frac{\left({F(R^2)-2\dot{F}(R^2)R^2}\right)^2}{{R}^2}$$
where $\dot{F}$ means the derivative of $F$ with  respect to the variable $t$.}
\end{theorem}
PROOF.
By using Proposition 3  we have that the extrinsic curvature of the surfaces of  rotation satisfies the following equation: 
\begin{equation}
\label{eqCC}[F+2\varphi \dot{F}(-\varphi+u\varphi')][F\varphi''-2a^2\dot{F}(-\varphi+u\varphi')]=-\tilde{K}_Ea^4\varphi,
\end{equation}
where $a^2=1+\varphi'^2.$

The surface $S(0,R)$ can be parametrized using   $ \varphi(u) = \sqrt{R^2 - u^2}$. In this case,
$$ -\varphi^{2} + u\varphi\varphi'= -R^2,\,  a^2( -\varphi + u \varphi') = {-R^4}/{\varphi^{3}} \mbox{ and }
  \varphi''= {-R^2}/{\varphi^{3}}.$$ 
  Substituting these expressions in $(\ref{eqCC})$ we obtain
\begin{equation}
\label{eqCCC}\tilde{K}_E= \frac{\left({F(R^2)-2\dot{F}(R^2)R^2}\right)^2}{{R}^2}. 
\end{equation}
This concludes the proof of Theorem.
 \hfill $ \Box $

In the next result we show that given any positive constant $c_0$ we can find a radius $R$, such that there 
exist a complete surface of rotation
in ${\mathbb{R}}^3_g$, with extrinsic curvature equal to $c_0$, that is a sphere centered at origin with  radius $R$.

\begin{theorem}
  {Let $c_0$ be a positive constant, $F$ and $\dot{F}$ are bounded functions globally defined in $\mathbb{R}$, with 
 $\displaystyle{ {\lim_{t\rightarrow \infty} \frac{F(t)-2\dot{F}(t) t}{\sqrt t}=0}}. $
Then there exist complete surface of rotation in ${\mathbb{R}}^3_g$,
 with extrinsic curvature equal to $c_0$, such surface is a sphere centered at origin with  radius $R$, 
 given by the relation} 
$$c_{0} = \frac{\left({F(R^2)-2\dot{F}(R^2)R^2}\right)^2}{{R}^2}. $$
\end{theorem}
PROOF.
Consider  $w(t)=\left(\frac{F(t)-2\dot{F}(t)t}{\sqrt {t}}\right)^2$.  By hypotheses $F(t)$ and $\dot{F}(t)$ are bounded
func\-tions, it fol\-lows that $\displaystyle{\lim_{t\rightarrow 0} 
w(t)= \infty}$. 

 Since that 
$\displaystyle{{\lim_{t\rightarrow \infty} w(t)=0}}$, then given any arbitrary constant
 $c_{0} >0$ there exist $R>0,$
such that 
$\displaystyle{c_{0} = \frac{\left({F(R^2)-2\dot{F}(R^2)R^2}\right)^2}{{R}^2},}$ 
proving that there exist a sphere centered at origin with radius $R$ in  ${\mathbb{R}}^3_g$, with extrinsic curvature $c_{0}$.  
 \hfill $ \Box $

In the next result we show that the straight lines through the origin are geodesics
\begin{theorem}
  {The radial lines in ${\mathbb{R}}^3_g$ are geodesics.}
\end{theorem}
PROOF.
Let us show that the straight lines through the origin are geodesics.
Let $\alpha (u)= u {\bf{v_0}}= u(v_1,v_2,v_3)$ be a  parametrized curve, which trace is a radial straight line,
 without loss of generality, we can consider $\bf v_0$ an unitary vector with respect to the Euclidean metric,
 i.e., $||\bf{v_0}||=1$. 
Hence $\alpha ' (u)=\bf v_0$, and 
$$||\alpha'(u)||_g^2=||{\bf v_0}||_g^2=\frac{||{\bf v_0}||^2}{F^2(\alpha(u))}.$$

  Let $s(u)=\int\limits_{0}^{u}||\alpha'(u)||_g du$ be the arc lenght function of $\alpha$, and denote 
by $\beta(s)=\alpha\circ h(s)$, the repara\-metriza\-tion by the arc lenght, where $h=s^{-1}$ is the inverse function of $s$.

 The system of ordinary differential equations of geodesics is given by (\ref{Geodesicas2}), namely
 $$
         \displaystyle \frac{d^2 x_k}{ds^2}+\frac{2x_k\dot F}{F}\left(\sum_{i\neq k}\left(\frac{dx_i}{ds}\right)^2-\left(\frac{dx_k}{ds}\right)^2\right) -\frac{4\dot F}{F}\frac{dx_k}{ds}
    \sum_{i\neq k} x_i \frac{dx_i}{ds}=0,
   $$
where $ k=1,2,3$.  

Now substituting the component functions of the curve in the system (\ref{Geodesicas2}), where $\frac{d^2 x_k}{ds^2}=h''(s) v_k$, and observing that $h'(s)=F\circ \beta(s)$, then the system is equivalent to 
$$\displaystyle h''v_k+2h v_k \frac{\dot F (\beta(s))}{F(\beta(s))}(-(h' v_k)^2+\sum_{i\neq k}(h' v_i)^2)
-4v_k\frac{h h'^2\dot F (\beta(s))}{F(\beta(s))} \sum_{i\ne k}v_i^2=0,$$
since we have that the equation $h''-2hF\dot F=0$ is verified, we can conclude that the straight lines through the origin are geodesics.
 \hfill $ \Box $


\begin{theorem}
  {The space ${\mathbb{R}}^3_g$ has the following properties:}
\begin{itemize}
\item[A]: { $S(0,R)$ have zero extrinsic curvature if, and only if, $F(t)= \sqrt{t}$;}
\item[B]: { $S(0,R)$ in ${\mathbb{R}}^3_g$ are totally geodesics if, and only if, $F(t)= \sqrt{t}$;}
\item[C]: { $S(0,R)$ in ${\mathbb{R}}^3_g$ have constant gaussian curvature for all the conformal factor $F$. Moreover,
all the spheres have Gaussian curvature equal to 1 if, and only if, 
$F(t)=\sqrt{t}$;}
\item[D]: { Circles centered at origin are geodesics in ${\mathbb{R}}^3_g$ if, and only if, $F(t)=\sqrt{t}.$;}
\item[E]: { If $F(t)=\sqrt{t}$, then the function $f(x)= \frac{x}{<x,x>_{g_0}}$ is an isometry.}
\end{itemize}
\end{theorem}
PROOF.
{\bf proof of Item A:}
It follows by $(\ref{eqCCC})$ that the Euclidean spheres $S(0,R)$ in ${\mathbb{R}}^3_g$
have zero extrinsic curvature $\tilde{K}_E= 0$, $\forall \, R>0$ if, and only if, $ F(t)- 2t\dot{F}(t)= 0$. This equation  has  as 
solutions $F(t)= \lambda\sqrt{t}$, for any nonzero constant $\lambda$. Without loss of generality we can consider $\lambda=1$.

\vspace{0.25cm}
\noindent
{\bf proof of Item B:} It follows immediately from the fact that the second fundamental form vanishes everywhere.
  
\vspace{0.25cm}
\noindent
{\bf proof of Item C:}   Let  $X(u,v)=(x_1(u,v),x_2(u,v),x_3(u,v))$ be a sphere of radius $R$ in ${\mathbb{R}}^3_g$.
Then the induced metric is given by 
$$<X_u,X_v>_g=\frac{1}{F^2}<X_u,X_v>_{{g_0}},$$
 where $g_0$ is the Euclidean metric and $F(t)=F(x_1^2+x_2^2+x_3^2)=F(R^2)=cte$. 
 We have that the coefficients of
 the first fundamental form, with respect to the conformal metric, are given by
$$\tilde{E}={E}/{F^2},\,\, \tilde{G}={G}/{F^2} \mbox{ and } \tilde{F}=0.$$
In this case, the Gaussian curvature $\tilde{K}_g$ is given by
$$\begin{array}{lcl}
\tilde{K}_g&=& \frac{-1}{2\sqrt{\tilde{E}\tilde{G}}}\left\{\left(\frac{\tilde{E}_v}{\sqrt{\tilde{E}\tilde{G}}}\right)_v+\left(\frac{\tilde{G}_u}{\sqrt{\tilde{E}\tilde{G}}}\right)_u\right\}\\
&=&F^2\frac{-1}{2\sqrt{{E}{G}}}\left\{\left(\frac{{E}_v}{\sqrt{{E}{G}}}\right)_v+\left(\frac{{G}_u}{\sqrt{{E}{G}}}\right)_u\right\}\\
&=&F^2/{R^2}.
\end{array}$$
 
Then  $\tilde{K}_g=1$, if, and only if, $F(t)= \sqrt {t}, \,\, t\geq 0$.
 
\vspace{0.25cm}
\noindent
{\bf proof of Item D:}
 Since the rotations are isometries, we will only show, without lost of generality, that the circles 
centered at origin,  in the plane $x_3=0$, parametrized by arc length 
$\beta(s)=(R \cos(Fs/R), R \sin(Fs/R),0),$ where $ F=F(R^2)$, are geodesics. 

In fact, the circle $\beta(s)=(R\cos(Fs/R), R\sin(Fs/R),0)$ is parame\-trized by arc lenght, and the com\-ponent 
functions $x_i(s)$ satisfy

 $$\begin{array}{lcl}
x_1'(s)&=&-F\sin(Fs/R)\\
x_2'(s)&=&F \cos(Fs/R )\\
x_1''(s)&=&-F^2/R \cos(Fs/R)\\
x_2''(s)&=&-F^2/R \sin(Fs/R).
\end{array}$$

 Then by using the expressions of the Christoffel's symbols (\ref{SimbChristoffel2}) we get that 
$\beta$  is geodesic if, and only if, the following system is satisfied

\begin{equation}
\label{Eq14}
\begin{array}{lcl}
-{F^2} \cos \left({Fs}/{R}\right) \left({1}/{R} -2R{\dot{F}}/{F}\right)&=&0\\
{-F^2}\sin\left({Fs}/{R}\right)\left({1}/{R} -2R{\dot{F}}/{F}\right)&=&0.
\end{array}
\end{equation}

The system (\ref{Eq14}) is satisfied if, and only if, $1 -2R^2{\dot{F}(R^2)}/{F(R^2)}=0$. The solutions, up to a multiplicative constant $\lambda$, are given by $F(t) = \sqrt t$.

\vspace{0.25cm}
\noindent
{\bf proof of Item E:} It follows by straightforward calculations comparing the first fundamental forms.
 \hfill $ \Box $

In the next result, we will show that the Radial Model $({\mathbb{R}}^3, <\;,\;>_g)$, where $F(t) = \sqrt t$,  has beautiful properties. 
\begin{theorem}
 {The Radial Model $({\mathbb{R}}^3, <\;,\;>_g)$, where $F(t) = \sqrt t$, has the following properties:
 it is a complete Riemannian manifold, and has non negative sectional curvature given by
$$K\left({\partial}/{\partial x_i},{\partial}/{\partial x_j}\right)(x)= \frac{x^2_{k}}{ x^2_{i} +x^2_{j} + x^2_{k}},$$
where $1\leq i, j ,k\leq 3$ are distincts;}
\end{theorem}
PROOF.
We will show that the divergent curves of the space have infinity length.
Let $ \gamma(u)=(x_1(u),x_2(u),x_3(u))$ be an arbitrary divergent curve and 
$\gamma'(u)=(x_1'(u),x_2'(u),x_3'(u))$ its tangent vector field. By using spherical coordinates 
$$x_1(u)=r(u) \sin(\theta)\cos(\phi), \,\,x_2(u)=r(u) \sin(\theta)\sin(\phi), \,\, x_3(u)=r(u) \cos(\theta),$$  
therefore $x_1^2+x_2^2+x_3^2=r^2(u)$, and 
$$x_1'^2+x_2'^2+x_3'^2=[r'(u)]^2+r^2[\theta'(u)+\sin^2(u)[\theta'(u)]^2].$$

Since $F(x_1^2+x_2^2+x_3^2)=\sqrt{x_1^2+x_2^2+x_3^2}$,  the $g$-length of $\gamma$ is given by
 $$
 \begin{array}{lcl}
 L_{\gamma}&=&\int\limits_{u_0}^{u_1} ||\gamma'(u)||_{g}du=\int\limits_{u_0}^{u_1}
 \frac{||\gamma'(u)||_{g_0}}{{(F\circ\gamma(u))}}du\\
 &&
 \\&=&\int\limits_{u_0}^{u_1} \frac{\sqrt{x_1'^2+x_2'^2+x_3'^2}}{{ r(u)}}du
 \\
 &&
 \\&=&\int\limits_{u_0}^{u_1} \frac{\sqrt{[r'(u)]^2+r^2(\theta'(u)+\sin^2(u)[\theta'(u)]^2)}}{{ r(u)}}du
 \\
 &\geq &   \int\limits_{u_0}^{u_1} \frac{\sqrt{[r'(u)]^2}}{{r(u)}}du =  \int\limits_{u_0}^{u_1} \frac{|r'(u)|}{r(u)}du
\geq |\int\limits_{u_0}^{u_1} \frac{dr}{du} \frac{1}{r(u)} du|
 \\
 &\geq& \int\limits_{u_0}^{u_1} \frac{dr}{du} \frac{1}{r(u)} du
 =  \int\limits_{u_0}^{u_1} \frac{dr}{r}=  ln(u_1)-ln(u_0).
\end{array} $$

If $\gamma$ is a divergent curve going to the origin, then $\displaystyle \lim_{u_0\rightarrow 0}\gamma(u_0)=(0,0,0)$,
 without loss of generality, we can assume $u_1=1$ then the $g$-length of $\gamma$ is 
$$\displaystyle L_{\gamma}\geq\lim\limits_{u_0\rightarrow 0}  (ln(1) - ln (u_0))=+\infty.$$ 

On the other hand, without loss of generality, we assume $u_0=1$,   if $\gamma$ diverges to the infinity, then
 we have that $$\displaystyle L_{\gamma}\geq\lim\limits_{u_1\rightarrow + \infty}  (ln(u_1) - ln (1))=+\infty.$$

Hence $\gamma$ has infinity length. Therefore the Radial Model $({\mathbb{R}}^3, <\;,\;>_g)$, where $F(t) = \sqrt t$, is a complete space.

 By using (\ref{Ksec}), the sectional curvatures are given by 
$$K\left({\partial}/{\partial x_i},{\partial}/{\partial x_j}\right)(x)= \frac{x^2_{k}}{ x^2_{i} +x^2_{j} + x^2_{k}},$$
where $1\leq i, j ,k\leq 3$ are distints. This concludes the proof of Theorem.
 \hfill $ \Box $

\begin{coro}
 {In the Radial Model $({\mathbb{R}}^3, <\;,\;>_g)$, where $F(t) = \sqrt t$, all the spheres centered at origin $S(0,R)$ are minimal surfaces.}
\end{coro}
PROOF.
 It follows from Theorem 3.4, that the spheres $S(0,R)$ are totally geodesic in the Radial Model $({\mathbb{R}}^3, <\;,\;>_g)$, where $F(t) = \sqrt t$, in particular the spheres $S(0,R)$ are minimal surfaces.
 \hfill $ \Box $

\begin{coro}
  {The unique straight lines in the Radial Model $({\mathbb{R}}^3, <\;,\;>_g)$, where $F(t) = \sqrt t,$ that are geodesics are the radials ones.} 
\end{coro}
PROOF.
Given the initial conditions, $\alpha(s_0)=P_0\neq O$ and $\alpha'(s_0)=v_0$ (non-null vector), if $OP_0$ is
a multiple of the vector $v_0$, then the radial line defined by the origin $O$ and $P_0$ is the geodesic.
Suppose that the line satisfying the initial conditions, $\alpha(s_0)=P_0$ and $\alpha'(s_0)=v_0$,
where $OP_0$ is not a multiple of the vector $v_0$, is a geodesic, since that there exist a circle 
centered at origin that have the line as the tangent line, by the uniqueness of solutions of the
ordinary differential equations that characterizes the geodesic this give us a contradiction.
  \hfill $ \Box $

\begin{coro}
 { The Radial Model $({\mathbb{R}}^3, <\;,\;>_g)$, where $F(t) = \sqrt t$, is isometric to the  product space ${\mathbb{S}}^2\times {{\mathbb{R}}} $, under the identification 
 $\Psi: {\mathbb{R}}^3\setminus \{0\}\rightarrow {\mathbb{S}}^2\times {{\mathbb{R}}}$ given by }

 \begin{equation}
\label{Eq15}
\Psi(x)=\left({x}/{||x||},
  \log ||x|| \right).
  \end{equation}
\end{coro}
PROOF.
 It follows from Theorems 3.3, 3.4, and 3.5.
 \hfill $ \Box $

We will call $({\mathbb{R}}^3, <\;,\;>_g)$, where $F(t) = \sqrt t$ the conformal model or Radial Model of ${\mathbb{S}}^2\times {\mathbb{R}}$.   

\section{Classes of Weingarten Surfaces in $\mathbb{S}^2\times \mathbb{R}$}

In this section we study two  classes of Weingarten surfaces in the {\em Radial Model}, which satisfy $\tilde{K}_E+\tilde{H}^2-\tilde{K}=0 $ and $2\tilde{K}_E-\tilde{K}=0 $, where $\tilde{K}$ is the Gaussian curvature, $\tilde{H}$ is the mean curvature and $\tilde{K}_E$ is the extrinsic curvature. Moreover, by using the relation between the curvatures of the {\em Radial Model} and the curvatures with respect to the euclidean metric ([CPS]), we  prove that  first class of the Weingarten  surfaces in {\em Radial Model} corresponds, up  to isometries, to the minimal surfaces in $\mathbb{R}^3$, and second class corresponds to EDSGHW - surfaces in Euclidean space $ \mathbb{R} ^ 3$. Consequently these two classes of surfaces  has a Weierstrass type representation depending  on two holomorphic functions.

Let $S$ be a regular surface in the standard euclidean space $(\mathbb{R}^3, g_0)$, we can without lost of generality consider that $S$ is parametrized by $X(u,v)$, where $ u,v$ are the isothermal parameters, i.e., $E(u,v)=G(u,v)$ and $<X_u, X_v>_{g_0}=0$.  In this case we have the following result.

\begin{proposition}
Let $X$ be a parametrized surface, where $E=G, <X_u,X_v>_{g_0}=0$, and let $K$ be the Gaussian curvature of $S$ with respect to the euclidean metric $g_0$, and $\tilde{K}$ the Gaussian curvature of $S$ with respect to the radial metric $g_F$, where $F=h(t), t=<X,X>$. Then 
\begin{equation}
\label{eq10}
\tilde{K}=F^2(t) K+\frac1E\left\{ 4 E(h(t)h''(t)-h'(t)^2)(t-(<X,N>)^2) + 4 E h(t)h'(t)(1+H<X,N>) \right\}.
\end{equation}
\end{proposition}

PROOF.
 By using the Gauss equation for isothermal parameters we have that
\begin{equation}
\label{eq1}
K=\frac{-1}{2\sqrt{EG}}\left\{  \left(\frac{E_v}{\sqrt{EG}}\right)_v+\left(\frac{G_u}{\sqrt{EG}}\right)_u\right \}=\frac{-1}{2E}\left\{\left( \frac{{E}_v}{{E}} \right)_v+\left( \frac{{G}_u}{{E}} \right)_u \right\}.
\end{equation}

We observe that $S=X(U)$ as a parametrized surface in $M_3=(\mathbb{R}^3, g_F)$ is also isothermal, namely, the coefficients of the first fundamental form satisfy
\begin{equation}
\label{eq3}
\tilde{E}=\frac{E}{F^2}=\tilde{G}=\frac{G}{F^2} \mbox{ and } <X_u,X_v>_{g_F}=\frac{1}{F^2} <X_u,X_v>_{g_0}=0.
\end{equation}

Then

\begin{equation}
\label{eq2}
\tilde{K}=\frac{-1}{2\sqrt{\tilde{E}\tilde{G}}}\left\{  \left(\frac{\tilde{E}_v}{\sqrt{\tilde{E}\tilde{G}}}\right)_v+\left(\frac{\tilde{G}_u}{\sqrt{\tilde{E}\tilde{G}}}\right)_u\right \}=\frac{-1}{2\tilde{E}}\left\{\left( \frac{\tilde{E}_v}{\tilde{E}} \right)_v + \left( \frac{\tilde{G}_u}{\tilde{E}} \right)_u  \right\}.
\end{equation}

Since that $g_F=\frac{1}{F^2}g_0$, then $\tilde{E}=\frac{E}{F^2}=\tilde{G}=\frac{G}{F^2}$ and the Gauss curvature is given by

\begin{equation}
\label{eq4}
\tilde{K}= F^2K+\frac1E\left\{ F(F_{uu}+F_{vv})- (F_u^2+F_v^2) \right\}.
\end{equation}

Consider $F(x)=h(t)$, where $x=(x_1,x_2,x_3)=X(u,v)\in S$, and $t=<X,X>_{g_0}=x_1^2+x_2^2+x_3^2$, (the conformal factor is radial), and $h$ is a positive differentiable function. 
Then $F_u=2h' <X,X_u>_{g_0}$ and $F_v=2h' <X,X_v>_{g_0}$. The second partial derivatives of $F$ are given by 

$\begin{array}{lcl}
F_{uu}&=& 4h''(<X,X_u>)^2+2 h' [<X_u,X_u>+ <X,X_{uu}>],\\ 
F_{vv} &=& 4h''(<X,X_v>)^2+2 h' [<X_v,X_v>+ <X,X_{vv}>].
\end{array}$

Therefore

\begin{equation}
\label{eq6}
F_{uu}+F_{vv}= 4 h'' [\left(<X,X_u>\right)^2+\left(<X,X_v>\right)^2]+2h' [2E+<X, X_{uu}+X_{vv}>],
\end{equation}

and

\begin{equation}
\label{eq7}
F_u^2+F_v^2= 4 h'^2[\left(<X,X_u>\right)^2 + \left(<X,X_v>\right)^2].
\end{equation}

Taking the $g_0-$orthonormal base for $\mathbb{R}^3$ the $\displaystyle\{\frac{X_u}{\sqrt{E}},\frac{X_v}{\sqrt{G}}, N \},$
then we can express 
$$
\displaystyle X=<X,\frac{X_u}{\sqrt{E}}>\frac{X_u}{\sqrt{E}}+ <X,\frac{X_v}{\sqrt{G}}>\frac{X_v}{\sqrt{G}}+<X,N>N
$$.

Hence 

$$
\begin{array}{lcl}
<X,X>&=&(<X,\frac{X_u}{\sqrt{E}}>)^2+(<X,\frac{X_v}{\sqrt{E}}>)^2+(<X,N>)^2\\
&=&\frac{(<X,{X_u}>)^2+
(<X,{X_v}>)^2}{E}+(<X,N>)^2
\end{array}.
$$

 Thus
\begin{equation}
\label{eq8}
\left(<X,X_u>\right)^2+\left(<X,X_v>\right)^2=E[<X,X>-\left(<X,N>\right)^2].
\end{equation}

 Moreover, 

\begin{equation}
\label{eq92}
<X_{uu}+X_{vv}, X>=<X,N>2EH.
\end{equation}

By substituting (\ref{eq6}--\ref{eq92}) into (\ref{eq4}) we get
\\
$\begin{array}{lcl}
\tilde{K}&=& F^2K+\frac1E\left\{ h[4h''E(<X,X>-(<X,N>)^2)+2h'(2E + 2EH<X,N>)]+\right.\\
&&\\
&&-
\left.4h'^2E(<X,X>-(<X,N>)^2)\right\}
\\
&=& F^2 K +\frac1E\left\{ 4E(hh''-h'^2)(<X,X>-(<X,N>)^2) + 4Ehh'(1 + H<X,N>) \right\}.
\end{array}
$

 This concludes the proof of proposition.
 \hfill $ \Box $

\begin{proposition} 
Let $F(t)=\sqrt t$, $t=<X,X>$ be the conformal factor, then the extrinsic curvatures, the square of the Mean curvature and the Gaussian curvature are related by

\begin{equation}
\label{eq12}
\tilde{K}_E=<X,X> K  +2<X,N> H +\frac{(<X,N>)^2}{<X,X>},
\end{equation}
\begin{equation}
\label{eq13}
\tilde{H}^2= <X,X> H^2+ 2 <X,N> H + \frac{(<X,N>)^2}{<X,X>}
\end{equation}

\begin{equation}
\label{eq14}
\tilde{K}=<X,X> K + 2  <X,N>H + 2\frac{(<X,N>)^2}{<X,X>}
\end{equation}
\end{proposition}
and thus we have that

\begin{equation}
\label{eq14.1}
\tilde{K}_E+\tilde{H}^2-\tilde{K}=H[<X,X> H +2 <X,N>].
\end{equation}
PROOF.
 Supposing  that $h(t)=\sqrt{t}$, we have that  $h'(t)=\frac{1}{2\sqrt t}$, and $h''=-\frac14 t^{\frac{-3}{2}}$. We need to evaluate $hh''-h'^2=-\frac{1}{2t},$ and $hh'=\frac12$. Substituting these into equation (\ref{eq10}) we obtain
\begin{equation}
\label{eq11}
\begin{array}{lcl}
\tilde{K}&=&t K+\frac{-2}{t} (<X,X>-(<X,N>)^2) + 4(\frac12)(1+H<X,N>)\\
&=&tK+2\frac{(<X,N>)^2}{<X,X>}+2H<X,N>,
\end{array}
\end{equation}
since $t=<X,X>$. 

 This concludes the proof of proposition.
 \hfill $ \Box $

\begin{theorem}
Let $S=X(U)$ be a regular parametrized surface in the standard euclidean space $(\mathbb{R}^3, g_0)$ and in $M_3$, where $ F(t)=\sqrt t, \,\,\,t=<X,X>.$ Then 

$\tilde{K}_E+\tilde{H}^2-\tilde{K}=0 $ in $M_3$ if, and only if, up to isometries, $X$ is a minimal surface in $\mathbb{R}^3$.
\end{theorem} 

PROOF.
It follows directly from  (\ref{matrizwH}),  (\ref{eq14.1}) and Theorem 3.4 (item [E]).  \hfill $ \Box $

As an immediate consequence of the theorem we have

\begin{coro}
The Weingarten surfaces of rotation in $\mathbb{S}^2\times\mathbb{R}$ obtained in the Theorem 4.1 are, up to isometries, the cathenoids. 
\end{coro}

\begin{theorem}
Let $S=X(U)$ be a regular parametrized surface in the standard euclidean space $(\mathbb{R}^3, g_0)$ and in $M_3$, where $ F(t)=\sqrt t, \,\, t = <X,X>.$ Then 

$2\tilde{K}_E -\tilde{K}=0 $ in $M_3$ if, and only if,   $S$ is a EDSGHW- surface in the euclidean space  $\mathbb{R}^3$.
\end{theorem} 

PROOF.
It follows directly from (\ref{eq12})  and (\ref{eq14}).  \hfill $ \Box $


\bibliographystyle{elsarticle-num}

\end{document}